\documentclass[11pt]{amsart}
\usepackage{amsmath,amsthm}
\usepackage{paralist}
\usepackage{graphics} %% add this and next lines if pictures should be in esp format
\usepackage{epsfig} %For pictures: screened artwork should be set up with an 85 or 100 line screen
\usepackage{graphicx}
\usepackage{epstopdf}%This is to transfer .eps figure to .pdf figure; please compile your paper using PDFLeTex or PDFTeXify.
\usepackage[colorlinks=true]{hyperref}
\hypersetup{urlcolor=blue, citecolor=red}

\newtheorem{rem}{Remark}
\title[Equilibria For the Aggregation Equation]{Explicit Equilibrium Solutions For the Aggregation Equation with Power-Law Potentials}

\author[Jos\'e A. Carrillo and Yanghong Huang]{}


 \email{carrillo@imperial.ac.uk}
 \email{yanghong.huang@manchester.ac.uk}

\thanks{$^*$ Corresponding author: xxxx}

\begin{document}
\maketitle

% Enter the first author's name and address:
\centerline{\scshape Jos\'e A. Carrillo}
\medskip
{\footnotesize
    % please put the address of the first author
    \centerline{Department of Mathematics, Imperial College London}
    \centerline{London SW7 2AZ, United Kingdom}
} % Do not forget to end the {\footnotesize by the sign }

\medskip

\centerline{\scshape Yanghong Huang}
\medskip
{\footnotesize
    % please put the address of the second  and third author
    \centerline{School of Mathematics, The University of Manchester}
    \centerline{Manchester M13 9PL, United Kingdom}
}

\bigskip

% The name of the associate editor will be entered by an editorial staff
% "Communicated by the associate editor name" is not needed for special issue.
\centerline{(Communicated by the associate editor name)}

%The abstract of your paper
\begin{abstract}
Despite their wide presence in various models in the study of collective behaviors,
explicit swarming patterns are difficult to obtain. In this paper,
special stationary solutions of the aggregation equation with power-law kernels
are constructed by inverting Fredholm integral operators or by
employing certain integral identities. These solutions are expected to
be the global energy stable equilibria and to characterize the generic behaviors
of stationary solutions for more general interactions.
\end{abstract}

%%%%%%%%%%%%%%%%%%%%%%%%%%%%%%%%%%%%

\section{Introduction}
The mathematical analysis of equilibria and traveling wave type patterns in many body descriptions of collective behavior models is one of the basic scientific problems in the explanation of coherent structures in mathematical biology and technology. They naturally appear in swarming of animal species, cell movement by chemotaxis, granular media interaction and self-assembly of particles, see for instance \cite{MR3143990,MR2744704,MR2256869} and the references therein. Consensus in orientation patterns have been reported in the seminal papers~\cite{levine2000self,d2006self} in which they numerically show the asymptotic stability of compactly supported groups of particles moving in a fixed direction. The amazing feature of these flock patterns is that the consensus in orientation can be established based only on attraction and repulsion effects. More precisely, the movement of each particle is assumed to follow the equations of motion
\begin{equation}\label{2ndorder}
    \frac{d}{dt}x_i=v_i,\quad
    \frac{d}{dt} v_i = f(|v_i|)v_i -\nabla_x\left(\frac{1}{N} \sum_{j\neq i}\nabla K(|x_i-x_j|)\right), \quad i=1,\cdots,N\,,
\end{equation}
where $f$ is a function modeling the self-propulsion
at low speed and friction at high speed, and $K$ is the radial pairwise interaction
potential. As the total number of individuals is large, the system of differential
equations is difficult to analyse and usually a continuum description based on mean-field limits is
adopted, either at the kinetic level for the particle distribution function ~\cite{carrillo2009double,MR2744704}
or at the hydrodynamic level for the macroscopic density and velocities \cite{MR2369988,carrillo2009double}.
Different stationary patterns were observed for the system at the equilibrium speed $s_0$ (with $f(s_0)=0$), for instance coherent moving flocks and single or double rotating mills~\cite{levine2000self,d2006self,carrillo2009double,MR3090592,CKR}. The stability of these patterns at the discrete particle level was analysed in \cite{MR3325085,MR3215070,MR3158478}. At the continuum level, these patterns are characterized by searching for continuous probability densities or probability measures $\rho$ of particle locations such that the total force acting on each individual balances out. This is equivalent to finding probability densities or measures $\rho$ such that
\begin{equation}\label{eq:basic}
    \nabla K\ast \rho = 0 \quad\mbox{ on supp}(\rho) \,.
\end{equation}
Being the problem posed on the support of the unknown density $\rho$ implies that the equation \eqref{eq:basic} is highly nonlinear. In fact, characterizing the interaction potentials $K$ such that these profiles are continuous or regular in their support has recently triggered lots of attention. The richness of the qualitative properties of particular steady solutions depending on the potential $K$ was reported in \cite{PhysRevE.84.015203}. They numerically investigate the shape of the most stable solution for the first order dynamics associated to \eqref{2ndorder} given by
\begin{equation}\label{1storder}
    \frac{d}{dt}x_i=-\nabla_x\left(\frac{1}{N} \sum_{j\neq i}\nabla K(|x_i-x_j|)\right), \quad i=1,\cdots,N\,,
\end{equation}
and they characterize several bifurcations starting from the particular solution of the Delta ring, a solution concentrated on a circle. However, getting conclusions about the qualitative properties of the continuum problem \eqref{eq:basic} turned out to be  mathematically involved. The existence of explicit formulas for flock and mill patterns for particular potentials such as the Morse potentials, originally used in \cite{d2006self}, was possible due to the particular properties of associated differential operators \cite{MR2533628,MR2788924,MR3143997,MR3251743}. Let us point out that mill patterns are also characterized by a probability density or measure satisfying a similar equation to \eqref{eq:basic} giving the balance between attractive-repulsive and centrifugal forces of the form $\nabla K\ast \rho = - s_0\nabla \ln |x| \mbox{ in supp}(\rho)$.

The continuum evolution model associated to the particle system \eqref{1storder} via the mean-field limit \cite{MR3331178} is given by
\begin{equation}\label{eq:aggreq}
\rho_t = \nabla\cdot(\rho\nabla K*\rho),
\end{equation}
where $\rho(t,x)$ is the mass density function or measure of particles/individuals in space at time $t\geq 0$. This equation is known as the aggregation equation and it has been proposed for swarming modeling in~\cite{mogilner1999non} as the inertia-less continuous system associated to \eqref{2ndorder}. This model has been thoroughly studied with attractive and repulsive-attractive potentials, see \cite{MR3331178,MR3393315} and the references therein. Depending on the singularity of the potential at the origin, the solutions enjoy different regularity properties or propagate certain regularity in time (even sometimes uniformly in time). All these early results imply that the regularity of the stationary states depends on the singularity or the strength of the repulsion of the potential $K$ at the origin. In other words, the stationary states' regularity is deeply connected with the dissipation properties of the evolution \eqref{eq:aggreq}. In fact, this equation has a gradient flow structure related to probability measures \cite{MR2053570,MR1964483,MR2769217,MR3143991}, whose implication for our problem here is that stable stationary states of the dynamics of \eqref{eq:aggreq} should be among the local minimizers in suitable topology of the interaction energy functional
\begin{equation}\label{interenergy}
    \mathcal{E}[\rho] = \int \rho (K*\rho) \,dx\,,
\end{equation}
defined over all positive measures of mass $M_0$.
Actually, it is shown \cite{MR3067832} that the local minimizers $\rho$ of the interaction energy \eqref{interenergy} in suitable transport distance topology have to satisfy the Euler-Lagrange equations:
\begin{equation}\label{eq:equirho}
    K\ast \rho = E  \quad\mbox{ in supp}(\rho) \qquad \mbox{and} \qquad K\ast \rho \geq E \quad\mbox{ in } \mathbb{R}^d
\end{equation}
with $E\,M_0=\mathcal{E}[\rho]$. The previous conditions have to be understood $\rho$-a.e for the first condition and a.e. Lebesgue for the second one in case $\rho/M_0$ is just a mere probability measure instead of a density function. Using these Euler-Lagrange conditions, it was shown that the dimension of the support of local minimizers depends strongly on the singularity  of the potential at origin \cite{MR3067832};  that local minimizers can be bounded, compactly supported or continuous on the support if the singularity is strong enough \cite{CDM}. The last point was shown by using a hidden connection discovered in \cite{CDM} between the Euler-Lagrange conditions \eqref{eq:equirho} and the classical obstacle problem for elliptic operators. Finally, the existence of global minimizers of the interaction energy has been obtained just recently under quite sharp conditions \cite{MR3356997,MR3336989}. Therefore, we know the existence of solutions to the Euler-Lagrange equations \eqref{eq:equirho} for a quite large family of potentials including the Morse potentials under the non H-stability condition \cite{d2006self,MR3356997} and all power law potentials
\begin{equation}\label{eq:kernel}
 K(x) = \frac{|x|^a}{a}-\frac{|x|^b}{b}
\end{equation}
with $a>b>-d$. Let us remind the reader that the convention $\frac{|x|^0}{0}:=\log (x)$ is used. Because of the simpler topology, the one-dimensional case is  in general better understood, see \cite{MR2755500,MR2782822,MR3268056} and the references therein.

In this paper, we consider certain special exact stationary solutions in the sense that they satisfy the first condition of the Euler-Lagrange equations \eqref{eq:equirho} $K\ast \rho = E$ on the support of density $\rho$ for certain ranges of the power-law potentials. In our examples, we will find compactly supported, radially symmetric densities $\rho$ which are continuous inside its support. The power-law potentials might be considered a very restrictive case, however they are the typical examples obtained by asymptotics expansions and they should show the generic boundary behavior near the support of the stationary states. They also allow direct connections to more classical problems involving the Laplacian or its fractional counterparts. In fact, certain range of interactions are intimately connected to fractional diffusions \cite{MR3393315,MR3239623}. Finally, having explicit solutions allows for excellent test cases for numerical schemes \cite{MR3023729,MR3372289}.

In this paper, we will always look for radial symmetric solutions
associated with the power-law kernel~\eqref{eq:kernel}, which are
supported on a ball $B_R$ and are continuous inside. Explicit solutions 
are constructed in the range of parameters corresponding 
to $-d<b<a$, $b\leq 2$ when either $a$ or $b$ is an even integer.
These solutions are derived from the condition $K\ast \rho = E$ on the support of 
density $\rho$ with mass $M_0$, and are reasonable candidate
to be  global minimizers of the interaction energy. 
The paper is organized as follows. 
Some basic integral equations and their solutions in one
dimension are discussed in Section 2, but detailed derivations 
are deferred in in the Appendix. The one dimensional case with either $a=2$ or $b=2$
is studied thoroughly in Section 3, where some key observations allow a unified 
yet simpler approach in multidimensional case treated in Section 4.
We finally conclude discussing possible extensions, and comment on certain open problems and conjectures 
of the generic behavior of solutions to \eqref{eq:equirho}.

%%%%%%%%%%%%%%%%%%%%%%%%%%%%%%%%%%

\section{Fredholm integral equations and integral identities}
\label{sec:prem}
To derive the exact steady states governed by~\eqref{eq:equirho} with
either $a$ or $b$ being an even integer,  the key is to invert the Fredholm
integral operator of the form
$$
\mathcal{L}[\rho](x):=\int_{B_R}
|x-y|^p \rho(y)dy\,,\qquad p\in (-1,\infty),
$$
that is, to find the solution of  $\mathcal{L}[\rho] = f$.
In one dimension, singular integral equations with power-law kernel
like this have been studied extensively~\cite{MR1728075}, usually
in connection with the theory of complex variables. The first step in
solving $\mathcal{L}[\rho](x)=f(x)$ is to differentiate both sides
of the equation such that the kernel $|x-y|^p$ becomes weakly singular,
with exponent $p$ between $-1$ and $0$. Depending on
the precise value of $p$, we have the following two cases.

If $p \in (2k-1,2k)$ for some non-negative integer $k$,
upon taking derivatives $2k$ times, the integral equation becomes
\begin{equation}\label{eq:abseq}
    \int_{-R}^{R} |x-y|^{-\nu} \rho(y)dy = f(x)
\end{equation}
with $\nu = 2k-p \in (0,1)$. The solution is given by (see~\cite{MR1728075} for the derivation)
\begin{align}\label{eq:abseqsol}
     \rho(x)=&\,\frac{\sin\pi\nu}{2\pi}\frac{d}{dx}\int_{-R}^x
     \frac{f(y)}{(x-y)^{1-\nu}}\,dy\\
     &-\frac{\cos^2\frac{\pi \nu}{2}}{\pi^2}
     \left(R^2-x^2\right)^{\frac{\nu-1}{2}} \mbox{P.V.}\!\!
     \int_{-R}^R \frac{(R^2-y^2)^{\frac{1-\nu}{2}}}{y-x}
     \left\{\frac{d}{dy}\int_{-R}^y \frac{f(z)}{(y-z)^{1-\nu}}dz
\right\}dy,\nonumber
\end{align}
including a Cauchy principal value integral in the second term.
Explicit solutions can be obtained for special right hand side $f$
like polynomials. For example, if $f(x) = 1$,
\begin{equation}\label{eq:abssol1}
\rho(x) = \frac{\cos \frac{\pi
\nu}{2}}{\pi}(R^2-x^2)^{\frac{\nu-1}{2}},
\end{equation}
and if $f(x)=x^2$,
\begin{equation}\label{eq:abssolss}
    \rho(x) = -\frac{2\cos \frac{\pi\nu}{2}}{\nu(\nu+1)\pi}(R^2-
    x^2)^{\frac{\nu+1}{2}} +
    \frac{\cos \frac{\pi \nu}{2}}{\pi\nu} R^2 (R^2-x^2)^{\frac{\nu-1}{2}}.
\end{equation}
The detailed derivation, mainly involving the evaluation
of the principal value integrals, is given in Appendix~\ref{sec:1stfht}.

If $p \in (2k,2k+1)$ for some non-negative integer $k$, upon taking
derivatives $2k+1$ times, the integral equation becomes
\begin{equation}\label{eq:signeq}
    \int_{-R}^R (x-y)|x-y|^{-\nu-1} \rho(y)dy = f(x)
\end{equation}
with $\nu = 2k+1-p\in (0,1)$. The solution is given by (see~\cite{MR1728075}
for the derivation)
\begin{align}\label{eq:signeqsol}
     \rho(x)=&\,\frac{c\sin\frac{\pi\nu}{2}}{\pi} (R^2-x^2)^{\frac{\nu}{2}-1}+
     \frac{\sin\pi\nu}{2\pi}\frac{d}{dt}\int_{-R}^x
     \frac{f(y)}{(x-y)^{1-\nu}}\,dy \\
     &+\frac{\sin^2\frac{\pi \nu}{2}}{\pi^2}
     \left(R^2-x^2\right)^{\frac{\nu}{2}-1}
     \mbox{P.V.}\!\! \int_{-R}^R \!\frac{(R^2-y^2)^{1-\frac{\nu}{2}}}{y-x}
     \left\{\frac{d}{dy}\int_{-R}^y \frac{f(z)}{(y-z)^{1-\nu}}dz
\right\}dy.\nonumber
\end{align}
In addition to similar principal value integrals as in~\eqref{eq:abseqsol},
the solution consists of a non-trivial null space spanned by the
function $(R^2-x^2)^{\frac{\nu}{2}-1}$. Explicit solutions can also be obtained
when $f$ is a polynomial. For example,
if $f(x) = x$, then
\begin{equation}\label{eq:signsols}
\rho(x) = \frac{\sin \frac{\pi \nu}{2}}{\nu\pi}
 (R^2-x^2)^{\frac{\nu}{2}}+
c\frac{\sin\frac{\pi\nu}{2}}{\pi} (R^2-x^2)^{\nu/2-1}\,, \qquad c\in\mathbb{R}\,.
\end{equation}

Besides the appearance of solutions to~\eqref{eq:abseq} or~\eqref{eq:signeq}, certain
integral identities have also to be established when the
in order to eliminate free parameters in the final steady states. In many cases, the energy has to
be calculated, to select the unique stable state.  These special integrals turn out to be
identities intimately connected to the solutions given above. For example,
from the solution~\eqref{eq:abssol1} for $f(x)=1$,
we obtain the identity
\begin{equation}\label{eq:inteqid0}
\int_{-R}^R |x-y|^{-\nu}(R^2-y^2)^{\frac{\nu-1}{2}} dy
= \frac{\pi}{\cos \frac{\nu\pi}{2}},
\end{equation}
valid for $\nu \in (-1,1)$ and $x\in (-R,R)$. Other identities with
higher exponents can be obtained by successive integrations.
For example, integrating twice on both sides of~\eqref{eq:inteqid0} gives
\[
    \int_{-R}^R |x-y|^{2-\nu}(R^2-y^2)^{\frac{\nu-1}{2}} dy
=\frac{(2-\nu)(1-\nu)\pi}{2\cos \frac{\nu\pi}{2}}x^2
+\frac{(1-\nu)\pi}{2\cos \frac{\pi \nu}{2}}R^2,
\]
where the integration constant (by evaluating at $x=0$) is represented
as special integrals from Appendix~\ref{sec:SF}. Since the exponent of the kernel $|x-y|$ is usually fixed in the problem, it is more convenient to
change $\nu$ to $\nu+2$ in the previous identity to have the
same interaction kernel $|x-y|^{-\nu}$, that is,
\begin{equation}\label{eq:inteqid2}
\int_{-R}^R |x-y|^{-\nu}(R^2-y^2)^{\frac{\nu+1}{2}}dy
=-\frac{\nu(\nu+1)\pi}{2\cos \frac{\pi \nu}{2}}x^2
+\frac{(\nu+1)\pi}{2\cos \frac{\pi \nu}{2}}R^2,
\end{equation}
which is valid for $\nu \in (-3,-1)$. In fact, this identity holds for a wider range
of parameter $\nu \in (-3,1)$, precisely when the integral on the right hand side
of~\eqref{eq:inteqid2} is well-defined. As we will see, these identities are more useful than the general solutions formula~\eqref{eq:abseqsol} and~\eqref{eq:signeqsol}, because the former can provide a unified approach for the construction of the steady states, regardless of the range of the parameters. When $p$ is an integer, the solution  for the integral equation~\eqref{eq:abseq} has to be derived differently, but usually can be
established equivalently from the above results by taking the limit. For this reason,
the special cases of integer exponents are not discussed separately in the following sections.

%%%%%%%%%%%%%%%%%%%%%%%%%%%%%%%%%%%%%%%%%%%%%%%%%%%%%%%%%%%
%%%%%%%%%%%%%%%%%%%%%%%%%%%%%%%%%%%%%%%%%%%%%%%%%%%%%%%%%%%

\section{Exact steady states in one dimension}
\label{sec:1d}
\begin{figure}[htp]
\begin{center}
\includegraphics[totalheight=0.25\textheight]{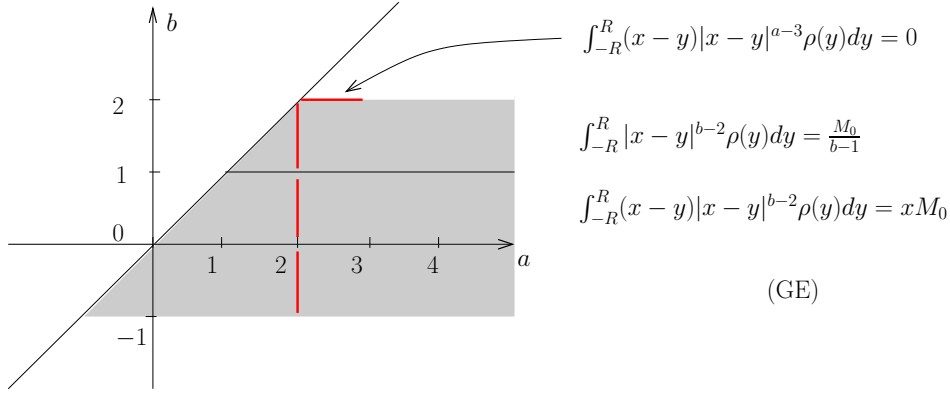}
\end{center}
\caption{The steady states exist only on the shaded
    region and their expressions are obtained from the reduced governing equation.}
\label{fig:1deq}
\end{figure}

In this section, we derive explicit expressions of the steady states in one dimension from the governing equation
\begin{equation}\label{eq:gov1d}\tag{GE}
\frac{1}{a}\int_{-R}^R |x-y|^a \rho(y)dy
-\frac{1}{b}\int_{-R}^R |x-y|^b \rho(y)dy = E, \qquad x\in (-R,R)\,,
\end{equation}
when either $a$ or $b$ is 2. These compactly supported symmetric steady states exist only in the
parameter regime $\{(a,b) | -1<b<2, b<a\}$ shown in Figure~\ref{fig:1deq}: the potential integral is well-defined only for $b>-1$, but the steady states becomes two delta masses as $b$ exceeds $2$ (see~\cite{MR3143991, MR2755500,MR2782822}). Upon repeated differentiation on both sides of~\eqref{eq:gov1d}, different integral equations of the form~\eqref{eq:abseq} or~\eqref{eq:signeq} appear, as summarized in Figure~\ref{fig:1deq}. In the special cases when either $a$ or $b$ is 2 of our interests, only the solutions~\eqref{eq:abssol1},~\eqref{eq:abssolss} and~\eqref{eq:signsols} are employed below.  The solutions are assumed to be radially symmetric and supported on some ball $B_R$ with unknown radius $R$. The precise expressions are derived first, and then the radius $R$ is determined from various self-consistency conditions like the definitions of certain moments. Without loss of general, we take the total mass $M_0$ given.

\subsection{The case when $b=2, a>2$}
We focus on the case $a \in (2,3)$ first. Taking derivative three times to~\eqref{eq:gov1d},
we get the homogeneous equation
\[
 \int_{-R}^R (x-y)|x-y|^{a-4}
\rho(y)dy =0\,.
\]
From~\eqref{eq:signeqsol}, the solution is given by
\begin{equation}\label{eq:solab2}
 \rho(x) = -c\frac{ \cos a\pi/2}{\pi}(R^2-x^2)^{\frac{1-a}{2}}.
\end{equation}
To determine the parameter $c$,  we first use the
definition of the total mass $M_0$, that is,
\begin{equation}\label{eq:massb2}
 M_0 = \int_{-R}^R \rho(x)\,dx
=- cR^{2-a}B\left(\frac{1}{2},\frac{3-a}{2}\right)
\frac{ \cos a\pi/2}{\pi}.
\end{equation}
On the other hand, certain information is lost during the differentiation and has to be retrieved upon substituting the
solution~\eqref{eq:solab2} back into~\eqref{eq:gov1d}.
By the special integral~\eqref{eq:symint} in Appendix~\ref{sec:SF}, we infer that
\[
\int_{-R}^R |x-y|^2 \rho(y)dy =
M_0x^2 + \int_{-R}^R y^2 \rho(y)dy =
M_0x^2 - c R^{4-a}B\left(\frac{3}{2},\frac{3-a}{2}\right)
\frac{ \cos a\pi/2}{\pi}.
\]
From the identity~\eqref{eq:inteqid2} with $\nu=-a$,
\begin{align*}
\int_{-R}^R |x-y|^a \rho(y)dy
= -c \frac{\cos a\pi/2}{\pi}
\int_{-R}^R |x-y|^{a} (R^2-y^2)^{\frac{1-a}{2}}ds
=\frac{c(a-1)}{2}\big( ax^2 + R^2).
\end{align*}
Therefore,  the left hand side of the governing equation~\eqref{eq:gov1d} is reduced to
\begin{equation}\label{eq:b2E}
 \left(\frac{c(a-1)}{2}-\frac{M_0}{2}\right)x^2 +
\frac{c(a-1)}{2a}R^2 +c R^{4-a}B\left(\frac{3}{2},\frac{3-a}{2}\right)
\frac{ \cos a\pi/2}{\pi}.
\end{equation}
For this expression to be a constant, the coefficient of $x^2$ should vanish, leading to the condition $M_0 = c(a-1)$.
This, together with~\eqref{eq:massb2},
can be simplified to determine the radius of the support
\begin{equation}\label{eq:1dRa2}
 R =\left[ -\frac{\cos a\pi/2}
     {\pi(a-1)}B\left(\frac{1}{2},\frac{3-a}{2}\right)
\right]^{\frac{1}{a-2}},
\end{equation}
and therefore the solution~\eqref{eq:solab2}.

\begin{figure}[htp]
\begin{center}
\includegraphics[totalheight=0.22\textheight]{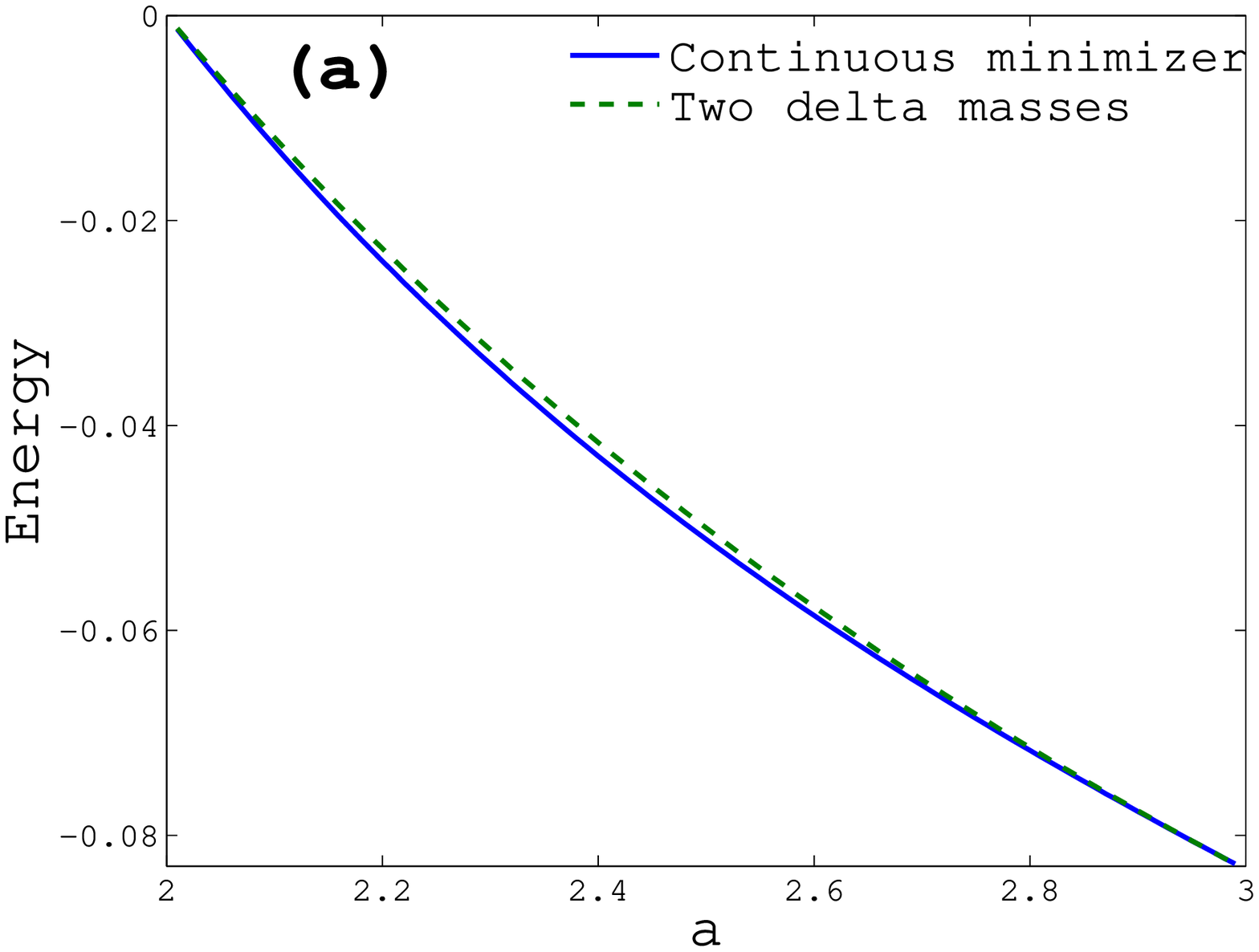}$~$
\includegraphics[totalheight=0.22\textheight]{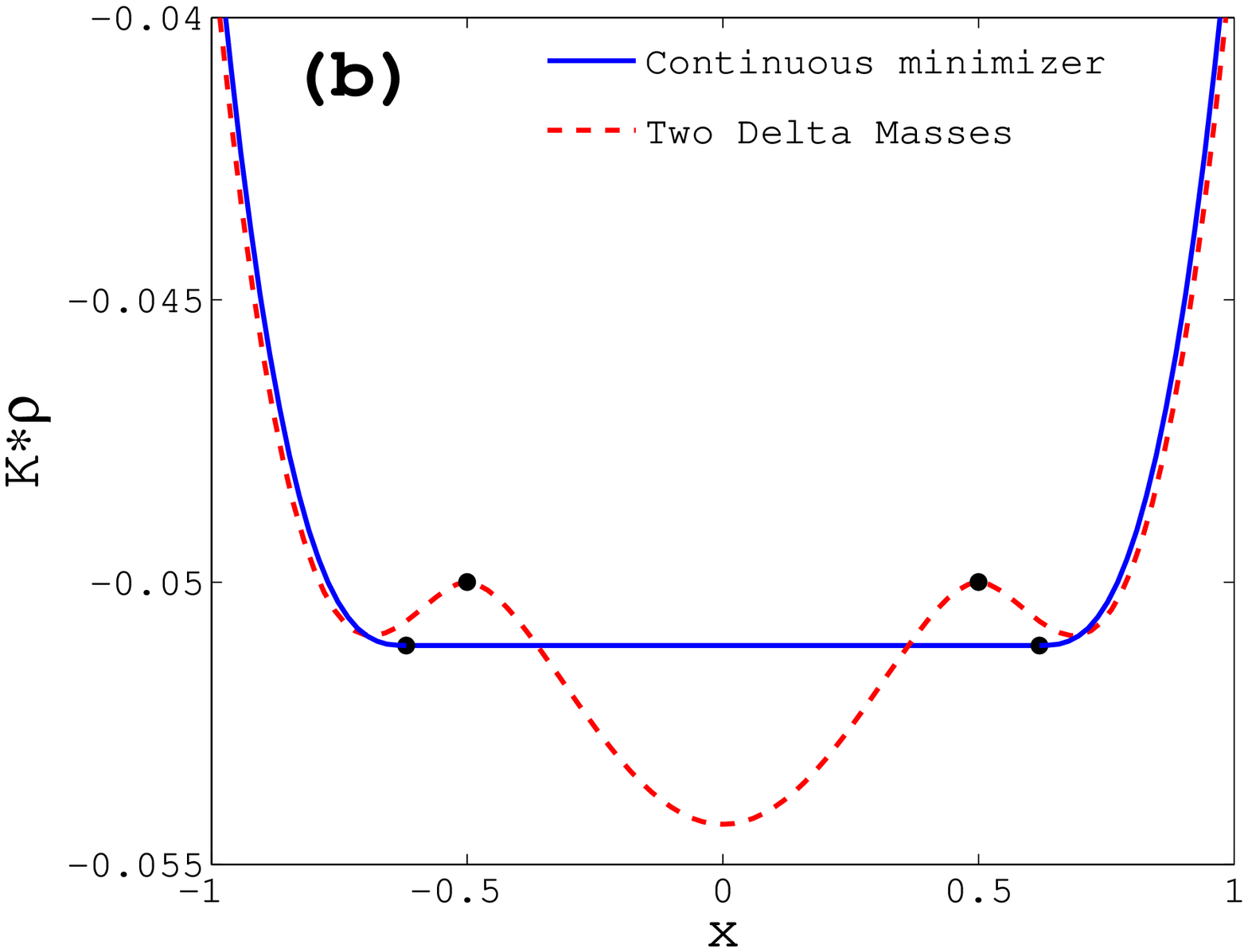}
\end{center}
\caption{{\bf (a)} the energy $E$ from the solution~\eqref{eq:solab2} compared with $E_\delta$ from two Dirac  masses; {\bf (b)} The function $K*\rho$ and $K\ast\rho_\delta$.}
\label{fig:compE}
\end{figure}

It is also instructive to calculate the associated energy in this case,
because the parameter  $b=2, a \in (2,3)$
lies on the boundary separating compactly supported solutions
that we are seeking and solutions consisting of two delta
masses~\cite{MR3067832,MR3325085}. From~\eqref{eq:b2E} and~\eqref{eq:1dRa2}, the energy
can be simplified as
\begin{equation}\label{eq:reEb2}
    E = \frac{c(a-1)}{2a}R^2 + c
    R^{4-a}B\left(\frac{3}{2},\frac{3-a}{2}\right)
    \frac{ \cos \frac{a\pi}{2}}{\pi}
    =\frac{2-a}{a(4-a)}M_0R^2.
\end{equation}
On the other hand, the more singular solutions consisting of two delta masses
can be written as
\begin{equation}\label{eq:delta1d}
    \rho_\delta(x) = \frac{M_0}{2}
    \Big( \delta_0(x-R_0) + \delta_0(x+R_0)\Big).
\end{equation}
Here the radius $R_0$ can be calculated in different ways, for example by minimizing the total energy
\[
  E_\delta:=\frac {1}{M_0} \int_{\mathbb{R}} \rho_\delta (K \ast \rho_\delta)\,dx
    =  \frac{M_0}{2} K(2R_0).
\]
This optimization procedure gives $R_0=1/2$ with the total energy $E_\delta =  (2-a)M_0/4a$. Though difficult to compare analytically, numerical
calculation shown in Figure~\ref{fig:compE}(a) implies that $E$ given by~\eqref{eq:reEb2} is always less than $E_\delta$ above, indicating the solution given by~\eqref{eq:solb2}
is more stable and it is believed to be the global stable solution in the sense of the global minimizer of the interaction energy \eqref{interenergy}, see \cite{MR3356997}.
The convolution $K\ast \rho$, shown in Figure~\ref{fig:compE}(b) for $a=5/2$, is constant
on $[-R,R]$ (the solid line between two dots) and is larger outside the support. In contrast, $K\ast \rho_\delta$ for $a=5/2$ behaves differently: although $R_0$ is the optimal radius,
the fact $K\ast \rho_\delta$ has a local maximum at $\pm R_0$ implies that the energy
can be reduced by making the density less concentrated.

We can proceed in the same way to find solutions for the other case when $a\geq 3$ and $b=2$, by differentiating the governing equation~\eqref{eq:gov1d}. If
$a \in (2k-1,2k)$ for some integer $k\geq 2$, the resulting homogeneous solution is of the form~\eqref{eq:abseq} and has only zero solution. If $a\in(2k,2k+1)$ for $k\geq 2$,
the solution $\rho(x)$ is proportional to $(R^2-x^2)^{(2k-1-a)/2}$ , whose coefficient
must be zero to make sure no higher order terms appear in~\eqref{eq:gov1d}.
In fact, densities like~\eqref{eq:delta1d} consisting two delta masses are the only stable steady states in this case.

\subsection{The case when $a=2, -1<b<2$}
For $a=2$, the governing equation becomes
\begin{equation}\label{eq:gov1da2}
 \frac{1}{2}\int_{-R}^R |x-y|^2 \rho(y) dy
-\frac{1}{b}\int_{-R}^R |x-y|^b \rho(y)dy = E.
\end{equation}
Depending on the value of $b$, derivatives of different orders
have to be taken.

If $b \in (1,2)$, we take derivative twice on~\eqref{eq:gov1da2}
to get
\begin{equation}\label{eq:1dgoveqb12}
 \int_{-R}^R |x-y|^{b-2} \rho(y)dy =\frac{M_0}{b-1}.
\end{equation}
From~\eqref{eq:abssol1} with $\nu=2-b$, the solution is given by
\begin{equation}\label{eq:solb2}
 \rho(x) = \frac{M_0}{b-1}\frac{\cos\frac{\pi(2-b)}{2}}{\pi}
(R^2-x^2)^{\frac{1-b}{2}}.
\end{equation}
Using the definition of the total mass
\[
 M_0=\int_{-R}^R \rho(x)dx
=\frac{M_0}{b-1}\frac{\cos\frac{\pi(2-b)}{2}}{\pi}
B\left(\frac{1}{2},\frac{3-b}{2}\right)R^{2-b},
\]
the radius of the support is then given by
\begin{equation}\label{eq:1dRb2}
 R = \left[
\frac{\cos \frac{(2-b)\pi}{2}}{(b-1)\pi}
B\left(\frac{1}{2},\frac{3-b}{2}\right)\right]^{\frac{1}{b-2}}.
\end{equation}

If $b \in (0,1)$, we take derivative once on~\eqref{eq:gov1da2} to get
\begin{equation}\label{eq:1dgoveqb01}
 \int_{-R}^R (x-y)|x-y|^{b-2}\rho(y)dy = xM_0.
\end{equation}
From~\eqref{eq:signsols} with $\nu=1-b$, the solution is can be written as
\begin{equation}\label{eq:1dsol01}
 \rho(x) = \frac{M_0\cos \frac{\pi b}{2}}{(1-b)\pi}
\Big[ (R^2-x^2)^{\frac{1-b}{2}}
    +cR^2(R^2-x^2)^{-\frac{1+b}{2}}
\Big]\,,
\end{equation}
where $c$ is the free parameter.
Because the component $(R^2-x^2)^{-\frac{1+b}{2}}$ is more
singular than $(R^2-x^2)^{\frac{1-b}{2}}$ near the boundary $|x|=R$, we must have the constraint $c\geq 0$ for the nonnegativity of the solution. The definition for the total mass is equivalent to
\begin{equation}\label{eq:a2b01}
 1 =  \frac{\cos \frac{\pi b}{2}}{(1-b)\pi}
 \left[ B\left(\frac{1}{2},\frac{3-b}{2}\right) +
c B\left(\frac{1}{2},\frac{1-b}{2}\right)
\right]R^{2-b},
\end{equation}
from which $c$ and $R$ are related.
Therefore, we have a one family of steady states parameterized by $c\geq 0$.
To select the most stable one of the family, we have to calculate the energy again
as a function of $c$.
The integration with the quadratic potential is straightforward,
\begin{align*}
\frac{1}{2}\int_{-R}^{R} |x-y|^2 \rho(y)dy &=
\frac{M_0}{2}x^2 + \frac{1}{2}\int_{-R}^R y^2 \rho(y)dy \cr
&=\frac{M_0}{2}x^2 + \frac{M_0\cos \frac{\pi b}{2}}{2(1-b)\pi}
\left[
    B\left(\frac{3}{2},\frac{3-b}{2}\right) +
    cB\left(\frac{3}{2},\frac{1-b}{2}\right)
    \right]R^{4-b}.
\end{align*}
Using~\eqref{eq:inteqid0} and~\eqref{eq:inteqid2} with $\nu=-b$, we get
\[
\int_{-R}^R |x-y|^b(R^2-y^2)^{-\frac{1+b}{2}}dy
=\frac{\pi}{\cos \frac{\pi b}{2}},
\]
and
\[
\int_{-R}^R |x-y|^b (R^2-y^2)^{\frac{1-b}{2}}dy
=\frac{b(1-b)\pi}{2\cos \frac{\pi b}{2}} x^2
+\frac{(1-b)\pi}{2\cos \frac{\pi b}{2}} R^2.
\]
Collecting all the terms above, the energy can be written as
\[
    E = \frac{M_0\cos \frac{\pi b}{2}}{2(1-b)\pi}
\left[ B\left(\frac{3}{2},\frac{3-b}{2}\right) +
cB\left(\frac{3}{2},\frac{1-b}{2}\right)
\right]R^{4-b}
- \frac{M_0R^2}{b(1-b)}\left[
c+\frac{1-b}{2} \right].
\]
Using the constraint~\eqref{eq:a2b01}, the energy can be simplified as
\begin{align*}
E &= \frac{M_0}{2}\frac{ B\left(\frac{3}{2},\frac{3-b}{2}\right) +
cB\left(\frac{3}{2},\frac{1-b}{2}\right) }{
    B\left(\frac{1}{2},\frac{3-b}{2}\right) +
    cB\left(\frac{1}{2},\frac{1-b}{2}\right) }R^2
-\frac{M_0R^2}{b(1-b)}\left[c+\frac{1-b}{2}\right]\cr
&= M_0\left[
\frac{1}{2}\frac{ (1-b)/(4-b)+c}{(1-b)+(2-b)c}
-\frac{c}{b(1-b)}-\frac{1}{2b}
    \right]R^2.
\end{align*}

Now, we show that $c=0$ is the global minimizer of $E$ for $c\geq 0$. First from~\eqref{eq:a2b01}, we can compute the derivative
\[
\frac{\partial R}{\partial c} = -\frac{R}{(2-b)c+1-b}.
\]
After some tedious algebra, we can show that
\[
\frac{d E}{d c} = \frac{\partial E}{\partial c} +
\frac{\partial E}{\partial R}\frac{\partial R}{\partial c}
=\frac{M_0(2-b)R^2c^2}{(1-b)(1-b+(2-b)c)^2}
>0
\]
for any $c >0$ and $b \in (0,1)$. Therefore, the energy $E$
is minimized precisely when $c=0$ and there is no singular component
$(R^2-t^2)^{-\frac{b+1}{2}}$ in the solution~\eqref{eq:1dsol01}.
The radius $R$ is given by the same formula~\eqref{eq:1dRb2}.

Finally, if $b\in (-1,0)$,
no derivative has to be taken and the original governing
equation \eqref{eq:gov1d} can be written as
\begin{equation*}
 \int_{-R}^R |x-y|^{b}\rho(y)dy
=\frac{M_0b}{2}x^2 + b\left(\frac{M_2}{2}R^2-E\right),
\end{equation*}
with the rescaled second order moment
$$
M_2=R^{-2}\int_{-R}^R |y|^2 \rho(y)\,dy\,.
$$
From~\eqref{eq:abssol1} and~\eqref{eq:abssolss} with $\nu=-b$,
the solution is given by
\begin{equation}\label{eq:1dsolb0-1}
    \rho(x) = \frac{M_0\cos \frac{\pi b}{2}}{(1-b)\pi}(R^2-x^2)^{\frac{1-b}{2}}
    +\left[\frac{bM_2-M_0}{2}R^2-bE\right]
    \frac{\cos \frac{\pi b}{2}}{\pi}(R^2-x^2)^{-\frac{1+b}{2}}.
\end{equation}
Define
\[
 c = \frac{1-b}{M_0R^2} \left[ \frac{bM_2-M_0}{2}R^2-bE\right],
\]
then by exactly the same procedure as for the case $b\in (0,1)$, one can show that the energy $E$ is an increasing function of $c$ for $c>0$. Therefore, the stable steady state is given by~\eqref{eq:solb2} for $c=0$ with the radius~\eqref{eq:1dRa2}.

\begin{rem}
From the constraint
\[
\frac{bM_2-M_0}{2}R^2-bE\geq 0
\]
for the coefficient of the second term in~\eqref{eq:1dsolb0-1}, it is tempting
to conclude that the energy is minimized if and only if this coefficient is zero.
However, when $E$ is taken as the only free parameter, both $R$ and $M_2$ solved from the definition of the first two even moments depend on $E$ too. As a result, it is not obvious that the
coefficient $(bM_2(E)-M_0)R(E)^2/2-bE$ is an increasing function of $E$. One illuminating
counterexample is to exam the case $b \in (0,1)$, where $-bE$ becomes an decreasing function of $E$ but the whole expression $(bM_2-M_0)R^2/2-bE$ has the opposite monotonicity.
\end{rem}

Although the approaches are different for $b$ in
different ranges of the interval $(-1,2)$, the final stable solution takes
the same form. This uniformity of the solution is a direct consequence of the
fundamental governing equation~\eqref{eq:gov1d}, since the integral equations~\eqref{eq:1dgoveqb12} and~\eqref{eq:1dgoveqb01} are obtained
from~\eqref{eq:gov1d} by taking derivatives. This observation is the key 
to derive general solutions in higher dimensions later.

\subsection{Extension with larger even integers of $a$ and in higher dimensions}

As we can see from the previous two subsections, the derivation of the
exact steady states is in general quite involved,  with different solutions
to integral equations and the establishment of integral identities in different parameter
regimes. The situation is further complicated by the appearance of free parameters that
usually have to selected by minimizing the energy.

The above approach using explicit solutions of the integral
equations~\eqref{eq:abseq} or~\eqref{eq:signeq}
cannot be easily generalized to higher dimensions either, because the
explicit dependence on angular integrals.  One exception is in three dimensions,
because of the following fact
\begin{align*}
\int_{B_R} |x-y|^p \rho(y) dy
&= 2\pi \int_0^{R} s^2\rho(s)\int_0^{\pi}
\big(r^2+s^2-2rs\cos\theta\big)^{p/2}\sin \theta d\theta ds \cr
&=\frac{2\pi}{(p+2)r}\int_0^R \big[ (r+s)^{p+2}-|r-s|^{p+2}\big]
s\rho(s)ds,
\end{align*}
with $r=|x|$ and $s=|y|$. This reduction is not useful yet, as
there are two nonlocal terms associated with one single power.
But these two terms can be combined into one integral by extending the radial
density $\rho(r)$ evenly on the interval $[-R,R]$, i.e.,
\[
    \int_{B_R} |x-y|^p \rho(y) dy
= -\frac{2\pi}{(p+2)r}\int_{-R}^R |r-s|^{p+2}s\rho(s)ds, \qquad s \in [-R,R].
\]
Therefore, the governing equation becomes
\begin{equation*}
-    \frac{1}{a(a+2)}\int_{-R}^R |r-s|^{a+2} s\rho(s)ds
+    \frac{1}{b(b+2)}\int_{-R}^R |r-s|^{b+2} s\rho(s)ds
=\frac{M_0E}{2\pi}r,
\end{equation*}
from which the steady states can be derived similarly as in one dimension, by taking
derivatives on both sides first. Since higher powers appear in the kernels, the derivative of the steady states is more cumbersome than that in one dimension, hence is not pursued here.

In summary, the above approach by solving integral equations like~\eqref{eq:abseq}
or~\eqref{eq:signeq} is getting more involved with larger exponents in one dimension, and
is unlikely to be generalized to higher dimensions.
However, certain features like the universal form of the solutions across different
parameter regimes indicate other unified ways, that we exploit next.

%%%%%%%%%%%%%%%%%%%%%%%%%%%%%%%%%%%%%%%%%%%%%%%%%%%%%%%%%%%%%%%%
%%%%%%%%%%%%%%%%%%%%%%%%%%%%%%%%%%%%%%%%%%%%%%%%%%%%%%%%%%%%%%%%

\section{Exact solutions in higher dimensions}
\label{sec:hd}

Now we focus on the construction of the  steady states in general dimensions,
based on some insights from the case of $a=2$ above: the solutions
take the same form in different parameter regimes and certain component never appears.

Assuming the radial solution $\rho$
is supported on the ball $B_R = \{ x\in \mathbb{R}^d, |x|\leq R\}$ and $a=2k$
is an even integer, the left hand side of the equivalent governing equation
\begin{equation}\label{eq:goveqhd0}
\int_{B_R}|x-y|^b \rho(y)dy = b\left[
\frac{1}{2k}\int_{B_R} |x-y|^{2k} \rho(y)dy-E
\right]
\end{equation}
is an even polynomial of degree $2k$, whose coefficients
can be written as moments of $\rho$.  It turns out that all
the solutions of~\eqref{eq:goveqhd0} can be built successively from the fundamental identity
\begin{equation}\label{eq:inteqidHD0}
    \int_{B_R} (R^2-|y|^2)^{-\frac{b+d}{2}}
    |x-y|^{b}dy = \frac{\pi^{\frac{d}{2}+1}}
    {\Gamma(\frac{d}{2})\sin\frac{(b+d)\pi}{2}}, \quad b\in (-d,2-d)
\end{equation}
which is well-known in potential theory~\cite[Appendix]{MR0350027}.
This identity (and the ones with higher exponents below) is expected from
dimensional analysis, by counting the powers of $R$ on both sides.

Once the identity~\eqref{eq:inteqidHD0} is established,
others with larger exponents can be derived by integration as those
in the end of Section~\ref{sec:prem}.
Denote $\Delta_x$ the Laplace operator, then
\begin{align*}
 \Delta_x  \int_{B_R} (R^2-|y|^2)^{-\frac{b+d}{2}}
    |x-y|^{b+2}dy
  &= (b+d)(b+2)\int_{B_R} (R^2-|y|^2)^{-\frac{b+d}{2}}
    |x-y|^{b}dy \cr
  & =     (b+d)(b+2)\frac{\pi^{\frac{d}{2}+1}}
    {\Gamma(\frac{d}{2})\sin\frac{(b+d)\pi}{2}}.
\end{align*}
By the radial symmetry of the integrals, the Laplace operator $\Delta_x$ can be
inverted easily, giving
\begin{equation*}
 \int_{B_R} (R^2-|y|^2)^{-\frac{b+d}{2}}
    |x-y|^{b+2}dy
= \frac{(b+d) \pi^{\frac{d}{2}+1}}  {\Gamma(\frac{d}{2})\sin\frac{(b+d)\pi}{2}}
\left( \frac{b+2}{2d}|x|^2 + \frac{1}{2}R^2\right),
\end{equation*}
where the integration constant can be obtained from special integrals
like~\eqref{eq:syminthd} in the Appendix, by evaluating both sides at the origin.
Replacing $b$ with $b-2$ to keep the same kernel $|x-y|^b$, we arrive at
\begin{equation*}
 \int_{B_R} (R^2-|y|^2)^{1-\frac{b+d}{2}}
    |x-y|^{b}dy
= -\frac{(b+d-2) \pi^{\frac{d}{2}+1}}  {\Gamma(\frac{d}{2})\sin\frac{(b+d)\pi}{2}}
\left( \frac{b}{2d}|x|^2 + \frac{1}{2}R^2\right)\,.
\end{equation*}
The range of validity for $b$ inherited from~\eqref{eq:inteqidHD0} is $(2-d,4-d)$,
but can be increased to $(-d,4-d)$ by continuation.
Integral identities with higher powers on $(R^2-|y|^2)$ can be
obtained similarly, which are special cases (for $k$ being integers) of the relation
\begin{multline}\label{eq:idhdk}
\int_{B_R} (R^2-|y|^2)^{k-\frac{b+d}{2}}|x-y|^b dy \\
=\frac{\pi^{\frac{d}{2}}}{\Gamma(\frac{d}{2})}
B\left(\frac{b+d}{2},k+1-\frac{b+d}{2}\right)
R^{2k}\, { }_2F_1\left(-\frac{b}{2},-k;\frac{d}{2};\frac{|x|^2}{R^2}\right).\quad 
\end{multline}
in view of the connection to the (negative) fractional Laplacian~\cite[Eq (9)]{MR3239623}.
When $k$ is a non-negative integer, the Gauss hypergeometric function ${ }_2F_1$ on the right hand side  is a
finite polynomial and~\eqref{eq:idhdk} can be written as
\begin{equation}\label{eq:inteqidHD}
    \int_{B_R} (R^2-|y|^2)^{k-\frac{b+d}{2}}|x-y|^b dy = C_{k0}R^{2k} + C_{k1}R^2|x|^{2k-2} + \cdots + C_{kk}|x|^{2k}
\end{equation}
with the following expressions
\begin{equation}\label{eq:IICoef}
 C_{kj} = \frac{\pi^{\frac{d}{2}}}{\Gamma(\frac{d}{2})}
    B\left(\frac{b+d}{2},k+1-\frac{b+d}{2}\right)
\frac{(-\frac{b}{2})_j(-k)_j}{j!(\frac{d}{2})_j}.
\end{equation}
Here the Pochhammer symbol $(r)_n=\Gamma(r+n)/\Gamma(r)=r(r+1)\cdots(r+n-1)$ is used.
The range of $b$ also increases with larger exponents, which
becomes $(-d,2+2k-d)$ for~\eqref{eq:inteqidHD}, exactly
when the components $(R^2-|y|^2)^{k-\frac{b+d}{2}}$ and $|x-y|^b$
of the integrand are integrable.

Equipped with these integral identities, we are now in a position to
construct the steady state for general even integer $a=2k$. First,
the governing equation~\eqref{eq:goveqhd0} can be written as
\begin{align*}
    \int_{B_R}|x-y|^{b}\rho(y)dy
    &= b\left[\frac{1}{2k}\int_{B_R}|x-y|^{2k}\rho(y)dy-E\right]\cr
    &= F_{0}R^{2k} + F_1R^{2k-2}|x|^2+\cdots + F_{k}|x|^{2k},
\end{align*}
for some coefficients $F_0,F_1,\cdots, F_{k}$ depending on the rescaled moments
$M_{2j}=R^{-2j}\int_{B_R}|x|^{2j}\rho(x)dx$. More precisely,
the coefficient $F_{j}$ only depends on the moments $M_0,M_2,\cdots, M_{2(k-j)}$,
but not on any higher order moments. From the integral identities above,
the steady states take the form
\begin{equation}\label{eq:rho}
\rho(x) =
(R^2-|x^2)^{-\frac{b+d}{2}}\Big(A_0R^{2k} + A_1R^{2k-2}(R^2-|x|^2)+\cdots
 + A_k(R^2-|x|^2)^k\Big),
\end{equation}
for some constants $A_0,A_1,\cdots,A_k$ to be determined. Here
appropriate powers of $R$ are used to make the coefficients $F_j$ and $A_j$
dimensionless and to simplify the calculations later.

\begin{figure}[htp]
\begin{center}
    \includegraphics[totalheight=0.25\textheight]{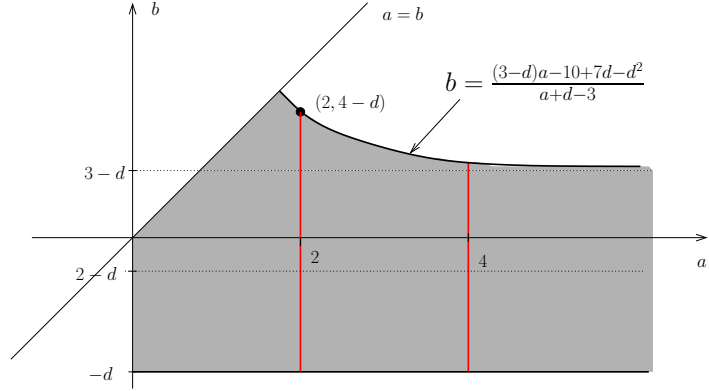}
\end{center}
\caption{The compactly supported steady states exist only in the shaded region.}
\end{figure}

Motivated by the solutions from the special case derived in the last section, we
now make some formal argument for $A_0$ vanishing identically in~\eqref{eq:rho}. First,
recall that we expect the radial steady states of our interest here to exist only when
\[
b<b_{max}:=\frac{(3-d)a-10+7d-d^2}{a+d-3},
\]
as determined in~\cite{MR3143991}. In fact, for $b>b_{max}$ the preferred asymptotic radial stationary state seems to be the uniform measure or Dirac Delta on a particular $(d-1)$-dimensional sphere \cite{PhysRevE.84.015203}, although it has not been rigorously proved.
Under the constraint $a>b$, it is easy to see that $b_{max}$ is always less than two and is exactly the reason
we only focus on $a=2k$, an even integer. Moreover, we have $3-d < b_{max} \leq 4-d$, with $b_{max}= 4-d$ only when $a=2$.

Now, we divide the range of $b$ into two intervals, $(2-d,b_{max})$ and $(-d,2-d)$.
For $b$ in the first interval $(2-d,b_{max}) \subset (2-d,4-d)$, the function $(R^2-|x|^2)^{-\frac{b+d}{2}}$
is no longer integrable near the boundary $|x|=R$. Therefore, $A_0$ must vanish identically
for the density to be well-defined. For $b$ belonging to the second interval
$(-d,2-d)$, the density $\rho$ is shown to be zero in~\cite{CDM} on the boundary $|x|=R$
by formulating the governing equation \eqref{eq:goveqhd0} as an equivalent
obstacle problem, which implies the coefficient $A_0$ of the unbounded yet integrable
component $(R^2-|x|^2)^{-\frac{b+d}{2}}$ must be zero.
In this range, the corresponding interaction energy is minimized by making a smooth transition to the zero
density outside the support, instead of creating an unbounded jump when $A_0$ is positive.

With this new information $A_0\equiv 0$, we can now construct the steady states
in a straightforward way. We show the calculation for $a=2$ and $a=4$ first
to demonstrate the procedure, and then go to general cases when $a$ is an even integer. For $a=2$, the governing equation is
\begin{equation}\label{eq:hdeq2}
\int_{B_R}|x-y|^b \rho(y) =b\left(
\frac{M_0}{2}|x|^2+\frac{M_2}{2}R^2-E
\right)
\end{equation}
The parameter $A_1$ in the solution $\rho(x) =
A_1(R^2-|x|^2)^{1-\frac{b+d}{2}}$ is determined uniquely
by matching the coefficient of $|x|^2$ on both sides of~\eqref{eq:hdeq2}, to obtain
\begin{equation}\label{eq:a2sol}
\rho(x) = -\frac{dM_0\Gamma(\frac{d}{2})\sin\frac{(b+d)\pi}{2}}
{(b+d-2)\pi^{\frac{d}{2}+1}}(R^2-|x|^2)^{1-\frac{b+d}{2}}.
\end{equation}
The radius $R$ is fixed from the definition of the total mass, that is
\[
    M_0 = \int_{B_R}\rho(x)dx = -\frac{M_0 d\sin\frac{(b+d)\pi}{2}}{(b+d-2)\pi}
    B\left(\frac{d}{2},2-\frac{b+d}{2}\right)R^{2-b}.
\]
In one dimension,~\eqref{eq:a2sol} reduces to the solution~\eqref{eq:solb2}.

When $a=4$, the governing equation becomes
\begin{equation}\label{eq:govhda4}
\int_{B_R} |x-y|^b\rho(y)dy
=F_0 + F_1R^2 x^2 + F_2|x|^4
\end{equation}
with
\[
F_0 = b\left(\frac{M_4}{4}R^4-E\right), \quad
F_1 = \frac{b(d+2)}{2d}M_2,\quad
F_2 = \frac{bM_0}{4}.
\]
Upon substituting
$\rho(x)=(R^2-|x|^2)^{1-\frac{b+d}{2}}(A_1R^2+A_2(R^2-|x|^2))$ followed
by the general identity~\eqref{eq:inteqidHD},
the matching condition of the coefficients of $|x|^2$ and $|x|^4$ on both sides of~\eqref{eq:govhda4}
becomes
\[
    \begin{pmatrix} C_{11} & C_{21} \cr
        0 & C_{22}
    \end{pmatrix}
    \begin{pmatrix} A_1 \cr A_2 \end{pmatrix}
        = \begin{pmatrix} F_1 \cr  F_2 \end{pmatrix}.
\]
Exacting $C_{kj}$ from~\eqref{eq:IICoef},  the solution of this linear system is given by
\begin{align}\label{eq:Asol}
A_1 &= \frac{\Gamma(\frac{d}{2})}{\pi^{\frac{d}{2}}}
\frac{1}{B(\frac{b+d}{2},2-\frac{b+d}{2})}\left[
\frac{d+2}{2}M_2 - \frac{d(d+2)}{2(b-2)}M_0
\right],\\
A_2 &= \frac{\Gamma(\frac{d}{2})}{\pi^{\frac{d}{2}}}
\frac{1}{B(\frac{b+d}{2},3-\frac{b+d}{2})}
\frac{d(d+2)}{4(b-2)}M_0,\nonumber
\end{align}
independent of $M_4$ and $R$. To decide the radius of the support, we use the
definition of the first two moments,
\[
M_0 = \int_{B_R} \rho(x)dx
= R^{4-b}\frac{\pi^{\frac{d}{2}}}{\Gamma(\frac{d}{2})}
\left[
    A_1 B\left(\frac{d}{2},2-\frac{b+d}{2}\right)
+A_2 B\left(\frac{d}{2},3-\frac{b+d}{2}\right)
\right]
\]
and
\begin{align*}
M_2 \,&= \frac{1}{R^2}\int_{B_R} |x|^2\rho(x)dx\\
&=R^{4-b} \frac{\pi^{\frac{d}{2}}}{\Gamma(\frac{d}{2})}
\left[
        A_1 B\left(1+\frac{d}{2},2-\frac{b+d}{2}\right)
        +A_2 B\left(1+\frac{d}{2},3-\frac{b+d}{2}\right)
        \right].
\end{align*}
From the explicit expressions~\eqref{eq:Asol} of $A_1$ and $A_2$, these two definitions
of the moments are homogeneous equations in $M_0$ and $M_2$, and
is equivalent to the eigenvalue problem
\begin{align*}
R^{b-4}
\begin{pmatrix} M_0 \cr M_2 \end{pmatrix} = \frac{(d+2)\Gamma(\frac{d}{2})}{2\Gamma(\frac{b+d}{2})\Gamma(2-\frac{b}{2})}
\begin{pmatrix} \frac{d}{4-b} & 1 \cr \frac{d^2}{(2-b)(6-b)} & \frac{d}{4-b} \end{pmatrix}
\begin{pmatrix} M_0 \cr M_2 \end{pmatrix}.
\end{align*}
The matrix on the right hand side has all positive entries. By Perron-Frobenius
theorem~\cite{MR1298430}, there is only one positive eigenvector $\vec{M} =[M_0,M_2]^T$
corresponding to the unique positive eigenvalue. In this case,
the special eigenvalue can be obtained explicitly, and the resulting radius is given by
\[
 R = \left[ \frac{ d(d+2)\Gamma(\frac{d}{2})}{2\Gamma(\frac{b+d}{2})\Gamma(2-\frac{b}{2})}
\left(\frac{1}{4-b}+\frac{1}{\sqrt{(2-b)(6-b)}}\right)
\right]^{-\frac{1}{4-b}}.
\]
The corresponding eigenvector can also be parameterized by the total mass $M_0$, that is
$M_2 = dM_0/\sqrt{(2-b)(6-b)}$ and the coefficient $A_1$ can be simplified as
\[
 A_1 = \frac{\Gamma(\frac{d}{2})}{\pi^{\frac{d}{2}}}
\frac{1}{B(\frac{b+d}{2},2-\frac{b+d}{2})}
\frac{d(d+2)}{2}\left[\frac{1}{\sqrt{(2-b)(6-b)}}+\frac{1}{2-b}\right]M_0.
\]
Therefore, the solution $\rho$ is completely determined.

Now we can proceed with the general case $a=2k$, and the  steps are outlined below:
\begin{enumerate}
 \item[Step a)] The solution takes the following form
 \begin{multline}\label{eq:genrho}
 \rho(x) = (R^2-|x|^2)^{-\frac{b+d}{2}}\Big[
 A_1R^{2k-2}(R^2-|x|^2) \cr
 +A_2R^{2k-4}(R^2-|x|^2)^2 +\cdots+A_k(R^2-|x|^2)^{k} \Big].
\end{multline}
\item[Step b)] Write $\vec{A}=[A_1,A_2,\cdots,A_k]^T$ as a function of $\vec{M}=[M_0,M_2,\cdots,M_{2k-2}]^T$ by matching
the coefficients of $|x|^2,\cdots, |x|^{2k}$ on both sides of the governing equation.
\item[Step c)] The definitions of the (rescaled) moments
$$
M_j=R^{-j}\int_{B_R} |x|^{j}\rho(x)\,dx\qquad \mbox{for } j=0,2,\cdots,2k
$$
can be formulated as an eigenvalue problem $R^{b-2k}\vec{M}
=\mathbf{D}\vec{M}$. The radius $R$ is the $1/(b-2k)$-th power of the eigenvalue of $\mathbf{D}$ and
the eigenvector, or equivalently the moments,  can be parameterized by the total mass $M_0$.
\end{enumerate}

\begin{rem}
Explicit forms of the coefficients can be worked out using symbolic packages for the general case $a=2k$.
But the calculation is more and more involved, and numerical softwares using floating point
are preferred.
\end{rem}

\begin{rem} The $k$-by-$k$ matrix $\mathbf{D}$ is expected to have positive entries,
leading to a unique eigenvalue $R^{2k-b}$ and a positive eigenvector $\vec{M}$. However, the
positivity of the moments does not imply that of the solution on the support,  as we will see below.
\end{rem}

\begin{figure}[htp]
\begin{center}
\includegraphics[totalheight=0.20\textheight]{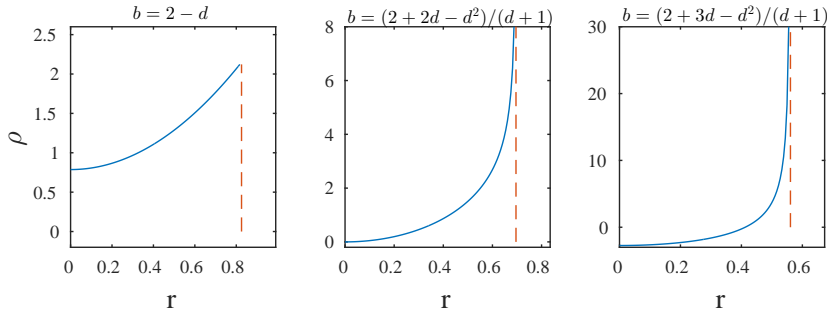}
\end{center}
\caption{The radial profiles of the constructed solutions for $a=4$. As $b$ increase from
 $2-d$ (the Newtonian potential) to $b_{max} = (2+3d-d^2)/(d+1)$, the solution
 becomes negative starting at $\bar{b}=(2+2d-d^2)/(d+1)$. }
 \label{fig:gapsol}
\end{figure}

However, the solutions constructed in this way are not always physically relevant, because
negative densities could appear. From the expressions~\eqref{eq:Asol} of $A_1$ and $A_2$ for $a=4$,
it is easy to see that the solution $\rho(x)=(R^2-|x|^2)^{1-\frac{b+d}{2}}(A_1R^2+A_2(R^2-|x|^2))$
becomes negative at the origin when
\[
 b > \bar{b} := \frac{2+2d-d^2}{d+1} \in (2-d,3-d).
\]
Since $\bar{b}$ is less than $b_{max}$, $\rho(0)$ is indeed negative for $b$ between $\bar{b}$ and $b_{max}$ (see Figure~\ref{fig:gapsol}). When $a$ is an even integer larger than $4$, numerical experiments show that
the solutions constructed in the above approach is still negative near the origin, when $b$ is close to its upper bound $b_{max}$. In fact, extensive particle simulations indicate that
a void (a region with zero density) starts to appear near the origin~\cite{MR3067832}.
As $b$ continues to increase, the density becomes concentrated more and more towards
its outer boundary, and eventually collapses on a sphere as $b$ passes $b_{max}$.
Because the solution is no longer supported on a ball, the above approach using integral identities like~\eqref{eq:idhdk} does not work, and precise form of the solutions remains unknown.

%%%%%%%%%%%%%%%%%%%%%%%%%%%%%%%%%%%%%%%%%%%%%

\section{Further discussions, conclusions and conjectures}

Besides providing more exact solutions for theoretical analysis and numerical 
testing,  these newly discovered exact solutions confirm many widely observed but yet 
unproved phenomena about boundary regularity even with the simplest power-law attractive-repulsive kernels.
From the general form~\eqref{eq:genrho}, it is easy to see that all constructed stationary solutions in this paper behave like $(R-|x|)^{1-(b+d)/2}$ near the boundary, depending only on the singular repulsive part of the kernel. More precisely, the Newtonian repulsive kernel dictated by $b=2-d$ is critical: the solutions is zero on its boundary for $b <2-d$ and become infinity (yet integrable) for $2>b>2-d$, while at the critical value $b=2-d$, the solution is finite on the boundary. We do expect this behavior to be generic for many potentials with power-law like repulsive component. As a result, if this boundary behavior conjecture is true, any converging solution of the evolution equation \eqref{eq:aggreq} can not be uniformed bounded in time for $2-d<b<2$.
Therefore, uniform (in time) $L^p$ bounds are expected only for $1\leq p < 2/(b+d-2)$.

The same boundary behaviors for other general potentials seem to be supported by extensive 
analytical and numerical studies. For Newtonian repulsion and general attraction, bounded density on the boundary is established as an eigenvalue problem  in~\cite{MR2834242,MR3143993}. Finite non-zero density on the boundary is proved for Newtonian-like repulsion like the Quasi-Morse potential~\cite{MR2788924}. 
The inverse square root singularity of the solutions for the Morse potential in two dimension (with $b=1$) is also confirmed numerically in~\cite{ryan2012analysis}. The boundary behaviors are rigorously proved for Newtonian or more singular than Newtonian repulsions and quite general attractions in \cite{CDM} using the equivalence to classical obstacle problems in elliptic equations of the necessary conditions \eqref{eq:equirho}.
This work may motivate further investigation towards the regularity of the solutions near boundary.

However, important questions like boundary behaviors are only inferred from special solutions and observed from 
limited numerical simulations. The extend to which these statements remain true are still widely open, 
and rigorous proofs may require sophisticated techniques from harmonic analysis and potential theory.
Because the solutions stop to be valid due to the appearance of negative density as $b$ is close to $b_{max}$,  the exact way in which the solution collapses to a co-dimension one sphere is not clear either. 
The construction of these nearly singular solutions is another interesting problem.

%%%%%%%%%%%%%%%%%%%%%%%%%%%%%%%%%%%%%%%%%%%%%%%%%%%%%%%%%%%%%%%%
%%%%%%%%%%%%%%%%%%%%%%%%%%%%%%%%%%%%%%%%%%%%%%%%%%%%%%%%%%%%%%%%
\appendix

\section{Special functions and special integrals}
\label{sec:SF}
Many special functions appear  during the derivation of the exact solutions, in particular the Euler
Gamma function
\begin{equation*}
    \Gamma(x)  = \int_0^{\infty} t^{x-1}e^{-t}dt,
\end{equation*}
and the Beta function
\begin{equation*}
    \mbox{B}(\alpha,\beta) = \int_0^1 t^{\alpha-1}(1-t)^{\beta-1}dt.
\end{equation*}
The Beta function can be represented using the Gamma function, that is $\mbox{B}(\alpha,\beta)
=\Gamma(\alpha)\Gamma(\beta)/\Gamma(\alpha+\beta)$. Certain
special integrals on the ball $B_R$ also appear repeatedly in this paper,
with
\begin{equation}\label{eq:symint}
\int_{-R}^R |x|^\alpha (R^2-|x|^2)^{\beta} dx
= R^{\alpha+2\beta+1}\,\mbox{B}\left(\frac{\alpha+1}{2},\beta+1\right),
\end{equation}
in one dimension and more generally
\begin{equation}\label{eq:syminthd}
\int_{B_R} |x|^\alpha (R^2-|x|^2)^\beta dx
=R^{\alpha+2\beta+d}\frac{\pi^{\frac{d}{2}}}{\Gamma(\frac{d}{2})}
B\left(\frac{\alpha+d}{2},\beta+1\right),
\end{equation}
in higher dimensions $\mathbb{R}^d$.

%%%%%%%%%%%%%%%%%%%%%%%%%%%%%%%%%%%%

\section{Derivation of the solutions to the singular integral equations}
\label{sec:1stfht}
When the right hand sides of two singular integral equations~\eqref{eq:abseq} and~\eqref{eq:signeq}
are polynomials, the solutions using the complicated formula~\eqref{eq:abseqsol} and~\eqref{eq:signeqsol}
can be simplified significantly. The key step is to evaluate the principal integral of the form
\begin{equation}\label{eq:FHT}
    \mbox{P.V.}\int_{-R}^R (R+y)^n\left(\frac{R-y}{R+y}\right)^{\frac{1-\nu}{2}}
    \frac{1}{y-x}dy
\end{equation}
for non-negative integers $n$ and $\nu\in (0,1)$. For example, if $f(x)=1$
in~\eqref{eq:abseq}, then
$$
\frac{d}{dx}\int_{-R}^x (x-y)^{\nu-1}f(y)dy
=(x+R)^{\nu-1}
$$
and the solution~\eqref{eq:abseqsol} is reduced to
\[
    \rho(x) = \frac{\sin \pi\nu}{2\pi} (x+R)^{\nu-1}
-\frac{\cos^2\frac{\pi\nu}{2}}{\pi^2}
(R^2-x^2)^{\frac{\nu-1}{2}}\mbox{P.V.}
\int_{-R}^R  \left(\frac{R-y}{R+y}\right)^{\frac{1-\nu}{2}}
\frac{1}{y-x}dy,
\]
with $n=0$ in~\eqref{eq:FHT}. Singular integrals like these
can be evaluated using \emph{Plemelj Formulae} for sectionally analytic
functions, which are analytic on the complex plane except a few
contours~\cite{MR1989049,MR1728075}.
\begin{figure}[htp]
    \begin{center}
    \includegraphics[totalheight=0.25\textheight]{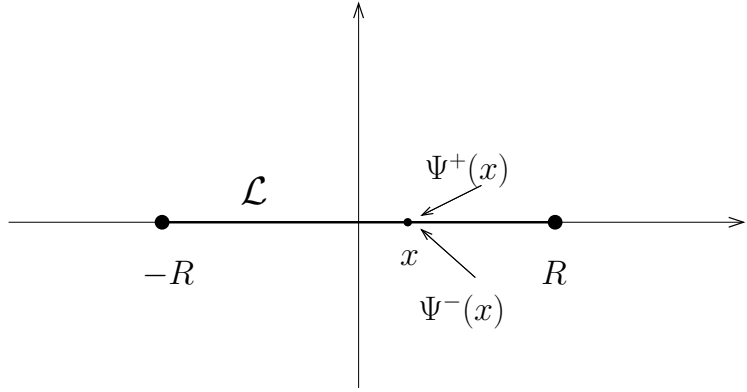}
    \end{center}
    \caption{The two limits of a sectionally analytical function $\Psi$ along
    the line segment connecting $-R$ and $R$.}
    \label{fig:plemelj}
\end{figure}

Let  $\mathcal{L}$ be the line segment connecting
$-R$ and $R$ and $\psi$ be a H\"{o}lder function on $\mathcal{L}$. Then the function defined by
\begin{equation}\label{eq:Psi}
    \Psi(z) = \frac{1}{2\pi i}\int_{\mathcal{L}} \frac{\psi(y)}{y-z}dy
\end{equation}
is analytic  on $\mathbb{C}\setminus{\mathcal{L}}$. For any real number
$x \in (-R,R)$, $\Psi$ has two limiting values when approaching $x$ from
above and below of $\mathcal{L}$, denoted as $\Psi^+(x)$ and $\Psi^-(x)$
(see Figure~\ref{fig:plemelj}).
The \emph{Plemelj Formulae} relates $\psi$ to the jump and average of
$\Psi^+(x)$ and $\Psi^-(x)$ in the following way,
\begin{subequations}
\begin{align}
    \psi(x)&=\Psi^+(x)-\Psi^-(x), \label{eq:plemejla}\\
    \frac{1}{\pi i}\mbox{P.V.}\int_{\mathcal{L}} \frac{\psi(y)}{y-x}dy
    &=\Psi^+(x)+\Psi^-(x). \label{eq:plemejlb}
\end{align}
\end{subequations}
Now we can calculate principal integrals like~\eqref{eq:FHT} using~\eqref{eq:plemejlb},
by finding appropriate complex function $\Psi$ whose jump across $\mathcal{L}$ is
exactly $(R+y)^n \left[(R-y)/(R+y)\right]^{(1-\nu)/2}$.

Let us choose $\Phi(z) =
\left(\frac{z-R}{z+R}\right)^{(1-\nu)/2}-1$ on $\mathbb{C}\setminus
{\mathcal{L}}$ with the branch cut $\mathcal{L}$. It is easy to calculate the limiting
values of $\Phi$ on $\mathcal{L}$ as
\[
    \Phi^\pm(x) = \left(\frac{R-x}{R+x}\right)^{(1-\nu)/2}
    e^{\pm i\pi(1-\nu)/2} -1\,.
\]
Finally, if we define
$$
\psi(x) := \Phi^+(x) - \Phi^-(x) =\left(\frac{R-x}{R+x}
\right)^{\frac{1-\nu}{2}} 2i\cos\frac{\nu\pi}{2},
$$
from the relation~\eqref{eq:plemejla}, the sectionally analytic function $\Psi$ given
by~\eqref{eq:Psi} has the same jump $\psi$ across $\mathcal{L}$ as $\Phi$.
Therefore, $\Psi-\Phi$ is continuous across $\mathcal{L}$ and  is analytic on the whole complex plane.
On the other hand, since both $\Psi$ and $\Phi$ are zero at infinity,  so is
their difference $\Psi-\Phi$ at infinity. Hence by Liouville's Theorem, $\Psi$ is identical to $\Phi$, and the
second Plemelj formula \eqref{eq:plemejlb} becomes
\[
\sin\frac{\nu\pi}{2} \left(\frac{R-x}{R+x}
\right)^{\frac{1-\nu}{2}}-1 =
\frac{\cos \frac{\nu\pi}{2}}{\pi}\,
\mbox{P.V.}\int_{-R}^R \left(\frac{R-y}{R+y}\right)^{\frac{1-\nu}{2}}
\frac{1}{y-x}dy.
\]
As a result, the principal value integral~\eqref{eq:FHT} with $n=0$ and $\nu\in (0,1)$ is obtained and the solution of~\eqref{eq:abseq} with $f(x)=1$ can be shown to
be~\eqref{eq:abssol1}. The same procedure can be generalized to larger
powers $n\geq 1$ and $\nu\in (0,1)$, by choosing
$$
\Phi(z)=(z+R)^n\left(\frac{z-R}{z+R}\right)^{\frac{1-\nu}{2}}-P_n(z),
$$
where $P_n(z)$ is a polynomial of degree $n$ such that $\Phi(z)$ becomes
zero at infinity. For example, we take
\[
\Phi(z) = (z+R)\left(\frac{z-R}{z+R}\right)^{\frac{1-\nu}{2}}
-z-\nu R
\]
to get (for $n=1$)
\[
\sin\frac{\nu\pi}{2} (R+x)\left(\frac{R-x}{R+x}
\right)^{\frac{1-\nu}{2}}-x-\nu R =
\frac{\cos \frac{\nu\pi}{2}}{\pi}
\int_{-R}^R (R+y)\left(\frac{R-y}{R+y}\right)^{\frac{1-\nu}{2}}
\frac{1}{y-x}dy\,,
\]
and
\[
\Phi(z) = (z+R)^2\left(\frac{z-R}{z+R}\right)^{\frac{1-\nu}{2}}
-z^2-(1+\nu) z R
-\frac{\nu^2+2\nu-1}{2}R^2
\]
to get (for $n=2$)
\begin{align*}
\sin\frac{\nu\pi}{2} (R+x)^2\left(\frac{R-x}{R+x}
\right)^{\frac{1-\nu}{2}}&\!\!\!\!-x^2 - (1+\nu) xR
-\frac{\nu^2+2\nu-1}{2}R^2 \\
&=
\frac{\cos \frac{\nu\pi}{2}}{\pi}
\int_{-R}^R (R+y)^2\left(\frac{R-y}{R+y}\right)^{\frac{1-\nu}{2}}
\frac{1}{y-x}dy\,.
\end{align*}
Once the principal value integral is obtained for any
non-negative integer $n$ and $\nu\in (0,1)$ in~\eqref{eq:FHT}, the
solution~\eqref{eq:abseqsol} can be evaluated for any polynomial left
hand side of~\eqref{eq:abseq}.

The solutions~\eqref{eq:signeqsol} for any polynomial right
hand side of~\eqref{eq:signeq} require the evaluation of a slightly different principal value integral
\begin{equation*}
    \mbox{P.V.} \int_{-R}^R (R+y)^n(R-y)\left(
    \frac{R+y}{R-y}\right)^{\frac{\nu}{2}}\frac{1}{y-x}dy,
\end{equation*}
for all non-negative integer $n$ and $\nu\in (0,1)$. This can be accomplished
by choosing the appropriate $\Phi$ in the Plemelj formula~\eqref{eq:plemejla} and~\eqref{eq:plemejlb}.
For example, we take
\[
 \Phi(z) = (z-R)^{1-\frac{\nu}{2}}(z+R)^{\frac{\nu}{2}}
-z - (\nu-1)R
\]
to get (for $n=1$)
\[
 -\cos \frac{\pi\nu}{2} (R-x)^{1-\frac{\nu}{2}}
(R+x)^{\frac{\nu}{2}} - x - (\nu-1)R =
\frac{\sin \frac{\pi\nu}{2}}{\pi}
\int_{-R}^R \frac{ (R-y)^{1-\frac{\nu}{2}}(R+y)^{\frac{\nu}{2}}}
{y-x}dy\,,
\]
and
\[
 \Phi(z) = (z-R)^{1-\frac{\nu}{2}}(z+R)^{1+\frac{\nu}{2}}
-z^2-\nu zR - \frac{\nu^2-2}{2}R^2
\]
to get (for $n=2$)
\begin{align*}
 -\cos \frac{\pi\nu}{2} (R-x)^{1-\frac{\nu}{2}}
(R+x)^{1+\frac{\nu}{2}}& - x^2 - \nu xR - \frac{\nu^2-2}{2}R^2 \nonumber\\
&=\frac{\sin \frac{\pi\nu}{2}}{\pi}
\int_{-R}^R \frac{ (R-y)^{1-\frac{\nu}{2}}(R+y)^{1+\frac{\nu}{2}}}
{y-x}dy\,.
\end{align*}

%For acknowledgements section, please don't number the section, please begin it with \section*{Acknowledgements}
\section*{Acknowledgments}
J.A.C. was partially supported by the Royal Society via a Wolfson Research Merit Award. J.A.C. and Y.H. were partially supported by the EPSRC grant EP/K008404/1.

\bibliography{biopowersol}
\bibliographystyle{plain}

\medskip
% The data information below will be filled by AIMS editorial staff
Received xxxx 20xx; revised xxxx 20xx.
\medskip

\end{document}